\documentclass[11pt]{article}
\usepackage{amsmath,amsthm,amssymb}
\newtheorem{theorem}{Theorem}
\newtheorem{lemma}{Lemma}

\newtheorem{remark}{Remark}
\theoremstyle{definition}
\newtheorem{definition}{Definition}
\begin{document}
\title{Bounded operators  on  weighted  spaces of holomorphic functions on the polydisk}
\author{A. V. Harutyunyan \\
Department of Applied Mathematics, Yerevan State University\\ 1 Alex Manookian str., 0025, Yerevan, Armenia\\
anahit@ysu.am}

\maketitle

\begin{abstract}We consider the weighted $A^p(\omega)$ and  $B_p(\omega)$  spaces of holomorphic functions on the polydisk (in the case of $p>1$). 
 We prove some theorems about    the boundedness of Toeplitz operators  on weighted Besov spaces $B_p(\omega)$ and about    the boundedness of generalized  little Hankel  and Berezin- type operators on $A^p(\omega)$.
\end{abstract}

\section{Introduction and auxiliary constructions.}
Let $U^n = \{ z = ( z_1, \ldots, z_n) \in \mathbb C^n, |z_j| <1, 1
\leq j \leq n \}$ be the unit polydisk in the n-dimensional complex
plane $ \mathbb C^n$ and let $T^n = \{ z = (z_1, \ldots, z_n) \in
\mathbb C^n, |z_j|=1, 1 \leq j \leq n \}$ be its torus. We denote by
$H(U^n)$ the set of holomorphic functions on $U^n$, by $L^{
\infty}(U^n)$ the set of bounded measurable functions on $U^n$  and by $H^\infty (U^n)$ the set of bounded holomorphic functions in $U^n$.

Let $S$ be the class of all non-negative measurable functions $
\omega$ on $(0,1)$ for which there exist positive numbers $M_{
\omega}$, $q_{ \omega}$, $m_{\omega}$ ($m_{ \omega}, q_{ \omega} \in
(0,1))$ such that
$$ m_{ \omega} \leq \frac{ \omega( \lambda r)}{ \omega(r)} \leq M_{
\omega } $$ for all $r \in (0,1)$ and $ \lambda \in [q_{\omega},
1]$. 
Some properties of the functions in $S$ can be found in \cite{se}. We put
\[ \alpha_{ \omega} = \frac{ \log m_{ \omega}}{ \log q^{-1}_{ \omega}};
\ \ \ \ \ \ \ \ \ \ \ \beta_{ \omega} = \frac{ \log M_{ \omega}}{
\log q_{ \omega}^{-1}}. \] For example $ \omega \in S$ if $
\omega(t) = t^{ \alpha} $ where $ -1 < \alpha < \infty$.  We
always have
$$ t^{ \alpha_{\omega}} \leq \omega(t) \leq t^{- \beta_{
\omega}},\,\,0<t<1. $$
 Below, for convenience of notations, for $\xi
= ( \xi_1, \ldots, \xi_n)$, $z=(z_1, \ldots, z_n)$ we put
$$\omega(1-|z|) = \prod_{j=1}^n \omega_j(1-|z_j|), \ \ \ 1-|z| =
\prod_{j=1}^n(1-|z_j|), \ \ \ 1- \bar{\xi} z = \prod_{j=1}^n (1 -\bar{\xi}_j z_j). $$ Furthermore, for $m = (m_1, \ldots, m_n)$, we
put $m+1 = \prod_{j=1}^n(m_j+1)$.

\begin{definition} Let $ 1< p < \infty$.
We denote by $L_p( \omega)$ the set of all measurable functions on
$U^n$ for which
$$||f||^p_{L^p( \omega)} := \int_{U^n} |f(z)|^p 
\omega(1-|z|)dm_{2n}(z) < \infty, $$

where $dm_{2n}(z)$ is the $2n-$dimensional Lebesgue measure on
$U^n$.
\end{definition}
We set $A^p(\omega)=L^p(\omega)\cap H(U^n)$ and denote the restriction of 
$|| \cdot ||_{L^P( \omega)}$ to $A^p( \omega)$ by $|| \cdot ||_{A^p( \omega)}$. In particular, for $\alpha=(\alpha_1,...,\alpha_n)$ if  $\omega_j(t)=t^\alpha_j,\,\,\alpha_j>-1\,\,(1\leq j\leq n)$ we set 
$A^p(\omega)\equiv A^p(\alpha).$

Note that the spaces
$A^p(\omega)$  are $\omega -$ weighted generalisations
of $A^p(\alpha)$ considered by F.Shamoyan (see \cite{d1}, \cite{sh}).

  The following definition gives the notion of the
fractional derivative.
 \begin{definition} (1) For a holomorphic
function $f(z) = \sum_{(k) = (0)}^{ (\infty)} a_k z^k$, $z \in U^n$
and for $ \beta = ( \beta_1, \ldots, \beta_n)$, $ \beta_j > -1$, $ 1
\leq j \leq n$,
  we define the fractional derivative $D^{ \beta}$ as follows
\[ D^{ \beta}f(z) = \sum_{(k) = (0)}^{( \infty)} \prod_{j=1}^n \frac{
\Gamma( \beta_j+1+k_j)}{ \Gamma( \beta_j+1) \Gamma(k_j+1)} a_kz^k,
\] $k=(k_1, \ldots,k_n)$, $z \in U^n$, where $ \Gamma( \cdot)$ is
the Gamma function and $ \sum_{(k) = (0)}^{ ( \infty)} =
\sum_{k_1=0}^{ \infty} \ldots \sum_{k_n=0}^{ \infty}$.

(2) Let $D^{- \beta}$ be the inverse operator, i.e. $D^{- \beta} D^{
\beta} f(z) = f(z) $ for $z \in U^n$.
\end{definition}
We put $Df(z) = D^{ \beta}f(z) $ if $ \beta=(1, \ldots, 1)$.

It is not difficult to show, that
\begin{equation}\label{e0}
f(z)=\int_0^1Df(rz)dr
\end{equation}
Next we define the holomorphic Besov spaces on the polydisk (see
\cite{hl}).
\vspace{4mm} \\
\begin{definition} Let $ 1 \leq p < \infty$ and $f \in
H(U^n)$. The function is said to be in $B_p( \omega)$ if
\[ ||f||^p_{B_p( \omega)} = \int_{U^n} |Df(z)|^p \frac{
\omega(1-|z|)}{(1-|z|^2)^{2 - p}}dm_{2n}(z) < \infty .\]\end{definition}

 We start by defining the
Toeplitz operator on the spaces $H(U^n)$. Let $L^1(T^n)$ be the
class of all integrable functions on $T^n$

\begin{definition} The Toeplitz operator with symbol $h
\in L^1(T^n)$ is the integral
\[ T_h(f)(z):= \frac{1}{(2 \pi i)^n} \int_{T^n} \frac{f( \xi) h(
\xi)}{ \xi -z } d \xi=\]
\[=\frac{1}{(2 \pi i)^n} \int_{T^n} \frac{f( \xi_1,\ldots,\xi_n) h(
\xi_1,\ldots,\xi_n)}{ (\xi_1 -z_1)\ldots (\xi_n-z_n)} d \xi_1\ldots
d\xi_n, \ \ \ f \in H(U^n). \]
\end{definition}

 \begin{remark} The above Toeplitz operator $T_h$ can
be extended to functions $f \in B_p( \omega)$ as follows: at first
one can consider $T_h$ on some everywhere dense subset of $B_p(
\omega)$, for instance on the set of all polynomials where $T_h$ is
obviously well-defined. Then one can show that, if the operator
$T_h$ is bounded with respect to $|| \cdot ||_{B_p( \omega)} $ on
the set of polynomials then it has a unique extension to $B_p(
\omega)$ naturally denoted by $T_h$ again.
\end{remark}
Our aim is to describe
the symbols $h$, for which $T_h$ defines a bounded operator $B_p(
\omega) \rightarrow B_p( \omega)$. 
Some problems about the Toeplitz operators can be solved by means of Hankel operators and vice versa. In the classical theory of Hardy of holomorphic functions  on the unit disk there is only one type of Hankel operator. In the $A^p(\alpha)$ theory they are two:  little Hankel operators and big  Hankel operators.  The analogue of the Hankel operators of the Hardy theory here are little Hankel operators, which were investigated by many authors (see for example  \cite{ja, af, z1}). 
Let us define the little Hankel operators as follows: denote by $\overline A^p(\omega)$ the space of conjugate holomorphic functions on $A^p(\omega) .$ For the integrable function $f$ on $U^n$  we define the generalized little Hankel operator  with symbol $h\in L^\infty(U^n) $ by
$$
h^\alpha_g(f)(z)=\overline P_\alpha (fg)(z)=\int_{U^n}\frac{(1-|\zeta |^2)^\alpha}{(1-\zeta\overline z)^{\alpha +2}}f(\zeta)g(\zeta )dm_{2n}(\zeta ),
$$
$$
\alpha =(\alpha_1,\ldots ,\alpha_n),\,\alpha_j>-1,\,1\leq j\leq n.
$$

The Berezin transform is the analogue of the Poisson transform in the $A^p(\alpha)$ theory. It plays an important role especially in the study of Hankel and Toeplitz operators. In particular, some properties of those operators (for example, compactness, boundedness) can be proved by means of the Berezin transform (see \cite{z1, st} ). The Berezin-type operators,  on the other hand, are of independent interest.

Next we consider generalized Berezin-type operators on $A^p(\omega)$ if $p>1.$  It will be  shown  that some properties of Berezin-type operators of the one dimensional classical case also hold in our more general situation.  For the integrable function $f$  on $U^n$  and for $g\in L^\infty (U^n)$ we define the Berezin-type operator in the following way

$$B_g^\alpha f(z)=\frac{ (\alpha +1)}{\pi^n}(1-|z|^2)^{\alpha
  +2}\int_{U^n}\frac{(1-|\zeta |^2)^\alpha}{|1-z\overline
  \zeta|^{4+2\alpha}}f(\zeta)g(\zeta )dm_{2n}(\zeta).$$
In the case $\alpha =0,\,\,g\equiv 1$ $B_g^\alpha $ will be called the Berezin transform. 

We need the following lemmas.
\begin{lemma}\label{l1}Suppose $f\in H^\infty (U^n)$ and $k=(k_1,\ldots,k_n),\,\, k_j \in N,\,\,1\leq j\leq n,$ then 
\begin{eqnarray*}
\left|\frac{\partial^{k_1+\ldots+k_n}f(z)}{\partial z_1^{k_1}\ldots\partial z_n^{k_n}}\right|\leq \frac{Const}{\prod_{j=1}^{n}(1-|z_j|)^{k_j}}
\end{eqnarray*}
\end{lemma}
{\bf Proof.} Let $f\in H^\infty (U^n).$  We consider the polydisk $$\widetilde T^n=\widetilde T_1\ldots\widetilde T_n,\quad\widetilde T_j=\left\lbrace \zeta_j,\zeta_j=z_j+\eta_j(1-|z_j|)e^{i\theta_j}\right\rbrace \quad j=1,...,n$$

We take $\eta _j$ so that $\widetilde T_j$ lies in the unit disc. Using the Cauchy formula for $\widetilde T^n$ we get
$$\frac {\partial^{k_1+\ldots+k_n} f(z)}{\partial z_1^{k_1}\ldots\partial z_n^{k_n}}=\frac{1}{(2\pi i)^n}\int_{\widetilde T^n}\frac{f(\zeta_1,\ldots\zeta_n )
  d\zeta_1\ldots d\zeta_n }{\prod_{j=1}^n(\zeta_j -z_j)^{k_j}},$$
then
\begin{eqnarray*}\Biggl|\frac {\partial^{k_1+\ldots+k_n} f(z)}{\partial z_1^{k_1}\ldots\partial z_n^{k_n}}\Biggr|\leq\frac{1}{(2\pi )^n}\int_{\widetilde T^n}
\frac{|f(\zeta )|d\theta_1...d\theta_n}
{\prod_{j=1}^n\eta_j(1-|z_j|^2)^{k_j}}\\
\leq Const\int_{\widetilde T^n}
\frac{d\theta_1...d\theta_n}{\prod_{j=1}^n\eta_j(1-|z_j|^2)^{k_j}}=
\frac{Const}{(1-|z_j|)^{k_j}}
\end{eqnarray*}

\hfill
$\square$

\begin{lemma}\label{l2}Let  $n=1, $ $\omega\in S ,\,\,\,a+1-\beta_\omega >0,$ $b>1$ and $ b-a-2>\alpha_\omega .$  Then
$$\int_{U}\frac{(1-|w|^2)^a \omega (1-|w|^2)}{|1-z\overline w|^{b}}dm_2(w)\leq\frac{\omega (1-|z|^2)}{(1-|z|^2)^{b-a-2}}$$\end{lemma}
For proof see Lemma 1.6 \cite{ha}

\section{Toeplitz operators on $B_p( \omega)$.}

\begin{theorem}  Let $p > 1, \,\,\beta_{\omega_j}<0$, 
 $p\geq\alpha_{\omega_j}$,$\beta_{\omega_j}+\alpha_{\omega_j}<0, \,\,1\leq j\leq n$ and  $h \in H^\infty (U^n)$.
  Then $T_{\overline {h}}$ is bounded operator: $B_p( \omega) \rightarrow
B_p( \omega).$
\end{theorem} 
{\bf Proof.}   
We show that, if $h \in B_q( \omega^*)$, then
$T_{\bar{h}}(f) \in B_p( \omega)$. Using 
\[ f(z) = C(\alpha,\pi)\int_{U^n} \frac{(1- |\xi|^2)^{ \alpha}Df(\xi)P(
\bar{\xi},z)}{(1- \bar{\xi}z)^{ \alpha+1}} dm_{2n}(\xi) . \] where
\[ P( \bar{\xi},z) = (1-(1- \bar{\xi}z)^{ \alpha+1})/z , \qquad
\alpha=(\alpha_1, ..., \alpha_n), \, \alpha_i\in \mathbb N\,, \]
 we get
\begin{eqnarray*}
&&T_{\bar{h}}(f)(z)= \frac{1}{(2 \pi i)^n} \int_{U^n}
(1-|\xi|^2)^m Df( \xi) \int_{T^n} \frac{ \overline{h(t)} P(
\bar{t},\xi) dt
dm_{2n}(\xi)}{(1- \bar{\xi}t)^{m+1}(t-z)} =\\
& & = \frac{1}{(2 \pi i)^n} \int_{U^n} (1-|\xi|^2)^m Df( \xi)
\overline{ \int_{T^n} \frac{ h(t) P( t,\bar{\xi}) t^{m} dt}{(1- t
\bar{z})(t-\xi)^{m+1}}  dm_{2n}(\xi)}
\end{eqnarray*}
Without loss of generality, we assume that
$P(\bar{t},\xi)=\bar{t}^l\xi^l$. Then we get
\begin{eqnarray*}
&&T_{\bar{h}}(f)(z)  = \frac{1}{(2 \pi i)^n } \int_{U^n}
(1-|\xi|^2)^m Df(\xi) \int_{T^n} \frac{\overline{h(t)}t^l
\xi^l\bar{t}^{m}dt}{(t-
\xi)^{m+1}{(1- \bar{t} z)}} dm_{2n}(\xi)\\
& & =\frac{1}{(2 \pi i)^n m!} \int_{U^n} (1-|\xi|^2)^m
Df(\xi)\xi^l\frac{\partial^m}{\partial\xi^m}\overline{\biggl(\frac{h(\xi)
\xi^{m-l}}{1-\xi\bar z}\biggr)}dm_{2n}(\xi)
\end{eqnarray*}
We set $ h(\xi)\xi^{m-l} = \tilde{h}( \xi)$. It follows that $\tilde h\in H^\infty(U^n)$ (we can take $m \in
\mathbb N^n$.)

  Then
\[ T_{\bar{h}}(f)(z) = \frac{1}{(2 \pi i)^nm!}
\int_{U^n}(1-|\xi|^2)^mDf( \xi) \sum_{|k|=0}^m \xi^kC_m^k \frac{
\partial^k \tilde{h}(\xi) }{ \partial \xi_1^{k_1} \ldots
\partial \xi_n^{k_n}}\frac{dm_{2n}( \xi)}{(1- \bar{\xi}z)^{m-k+1}}. \]
Let
\[ \Phi_k(z) = \int_{U^n} \frac{(1-| \xi|^2)^mDf( \xi)\xi^k}{(1-
\bar{\xi}z)^{m-k}} \frac{
\partial^{|k|} \tilde{h}(\xi)}{\partial \xi_1^{k_1} \ldots \partial
\xi_n^{k_n}}dm_{2n}( \xi). \] We shall show that $\Phi_k(z) \in B_p(
\omega)$.  To this end
using Lemma \ref{l1} we have
\begin{eqnarray*}
&&|D\Phi_k(z)|   \leq \int_{U^n}
\frac{(1-|\xi|^2)^m}{|1-\bar{\xi}z|^{m-k+2}} \left| \frac{
\partial^{k_1+\ldots+k_n} \tilde{h}( \xi, \bar{\xi})}{\partial
\xi_1^{k_1} \ldots \partial \xi_n^{k_n}} \right||Df(\xi)| dm_{2n}(\xi) =\\
  & & =
\int_{U^n}
\frac{(1-|\xi|^2)^{m-k}(1-|\xi|^2)^k}{|1-\bar{\xi}z|^{m-k+2}} \left| \frac{ \partial^{|k|} \tilde{h}(
\xi, \bar{\xi})}{\partial \xi_1^{k_1} \ldots \partial \xi_n^{k_n}}
\right| Df(\xi)| \leq\\
  & & \leq C||
\tilde{h}||_{\infty} \left( \int_{U^n}
\frac{(1-|\xi|^2)^{(m-k+2/q)p}|Df(\xi)|^pdm_{2n}(\xi)}{|1-\bar{\xi}z|^{(m-k
+2)p}} \right)^{1/p}
\end{eqnarray*}
Then
\[ \int_{U^n} (1-|z|^2)^p \frac{ |D \Phi_k(z)|^p
\omega(1-|z|)}{(1-|z|^2)^2}dm_{2n}(z) \leq ||h||_{\infty}\times\]
\[\times\int_{U^n}(1-|\xi|^2)^{(m-k+2/q)p}|Df(\xi)|^p
 \int_{U^n} \frac{(1-|z|^2)^p\omega(1-|\xi|)
dm_{2n}(\xi)}{|1- \bar{\xi}z|^{(m-k+2)p}(1-|z|^2)^2} \]
We estimate the following integral in the one dimensional case
\[ I( \xi) = \int_{U^n} \frac{(1-|\xi|^2)^{p-2}
\omega(1-|\xi|)}{|1-\bar{\xi}z|^{(m-k+2)p}}dm_2(z) .
\]
(i) If $(m-k+2)p-1-p+2+ \beta_{\omega} < 1$  (hence $(m-k+1)p < -
\beta_{ \omega}$) then
$I(\xi) \leq  const$. \\
(ii) If $(m-k+2)p > - \beta_{\omega}$ then $I(\xi) \leq
const./(1-|\xi|)^{(m-k+1)p+ \beta_{ \omega}-1}$. \\
(iii)If $(m-k+2)p = - \beta_{ \omega}$ then $I(\xi) \leq
\log(1-|\xi|)^{-1}$.\\
For the case (i) we have
\begin{eqnarray*}&&\int_{U}(1-|\xi|^2)^{(m-k+2/q)p}|Df(\xi)|^p \frac{
\omega(1-|\xi|)}{\omega(1-|\xi|^2)}dm_2(\xi) \leq \\
&&\int_U |Df(\xi)|^p
\frac{ \omega(1-|\xi|)}{(1-|\xi|^2)^{2-p}}(1-|\xi|)^\gamma dm_2(\xi)=\|f\|^p_{B_p(\omega)}
\end{eqnarray*}
 where $\gamma=p(m-k)+2p/q+2-p-\alpha_\omega =(m-k)p+p-\alpha_\omega
 \geq 0$.\\
For the case (ii) we have
\begin{eqnarray*}
&&\int_{U}(1-|\xi|^2)^{(m-k+2/q)p}|Df(\xi)|^p 
\frac{dm_2(\xi)}{(1-|\xi|)^{(m-k+2)p+
\beta_{ \omega}-1}} =\\
  & & = \int_U |Df(\xi)|^p \frac{ \omega(1-| \xi|)}{(1-|\xi|^2)^{2-p}}
\frac{dm_2( \xi)}{(1-|\xi|^2)^{ \beta_{\omega}+2p-2p/q-2}} \leq
C||f||_{B_p( \omega)}^p,
\end{eqnarray*}
if $ \beta_{ \omega} + 2p - 2 p/q -2 = \beta_{ \omega}+\alpha_\omega < 0$.\\
(iii) let $(m-k+1)p = - \beta_{ \omega}$. Then
\begin{eqnarray*}
&& \int_{U}(1-|\xi|^2)^{(m-k+2/q)p}|Df(\xi)|^p \frac{
\omega(1-|\xi|)}{\omega(1-|\xi|^2) }\log \left( \frac{1}{1-| \xi|}
\right)
dm_2( \xi)  \\
  & & =  \int_U |Df(\xi)|^p 
  \frac{ \omega(1-|\xi|)}{(1-|\xi|^2)^{2-p}}(1-|\xi|^2)^{ -\beta_{\omega}-p+2p/q-p+2-\alpha_\omega}
\log\left( \frac{1}{1-|\xi|} \right) dm_2( \xi)\\
& &  \leq C||f||_{B_p( \omega)}^p,\,\,
 -\beta_{ \omega}-\alpha_{\omega} > 0.
\end{eqnarray*}

 Then it follows that $ \Phi_k \in B_p(
\omega)$ and therefore $T_{ \bar{h}}(f) \in B_p( \omega) $.

 \hfill
 $\square$

We turn to the application of our result to division theorem
in spaces $B_p(\omega),\,(p> 1)$. To this end, we need the
following well-known definitions.

 \begin{definition} A function $g\in H^{\infty}(U^n)$ is
called an inner function, if its radial boundary values satisfy
$|g^*(w)|=1$ almost everywhere on $T^n$ (see \cite{ru}
\end{definition}
\begin{definition}An inner function  $g\in
H^\infty(U^n)$ is said to be good, if $u[g]=0$, where $u[g]$ is the
least $n$-harmonic majorant of $\log|g|$ in $U^n$ .
\end{definition}
\begin{theorem}  Let $p > 1, \,\,\beta_{\omega_j}<0$, 
 $p\geq\alpha_{\omega_j}$,$\beta_{\omega_j}+\alpha_{\omega_j}<0, \,\,1\leq j\leq n.$ Let $f\in B_p(\omega)$ and $J$ be good
inner function and
  $F=f/g\in H(U^n).$ Then $F\in B_p(\omega)$.
\end{theorem}

 {\bf Proof } is ivident:
$$
T_{\overline J}(f)(z)=\frac{1}{(2\pi i)^n}\int_{T^n}\frac{f(\zeta
  )\overline J(\zeta)}{\zeta -z}d\zeta=$$
$$=\frac{1}{(2\pi i)^n}\int_{T^n}\frac{f(\zeta
  )/{\overline J(\zeta)}}{\zeta -z}d\zeta=\frac{1}{(2\pi i)^n}\int_{T^n}\frac{F(\zeta
  )}{\zeta -z}d\zeta=F(z)$$
By Theorem 1  $F\in B_p(\omega)$

 \hfill
 $\square$

\section{Hankel und Berezin type operators on $A^p(\omega)$}

\begin{theorem}\label{th}Let $p> 1,\,\, f\in A^p(\omega)$ (or $f \in\overline A^p(\omega)) ,\,g\in L^\infty (U^n).$  Then $h^\alpha_g(f)\in \overline A^p(\omega)$  if and only if $\alpha_j>(\alpha_{\omega_j}+1)/p-1,\,1\leq j\leq n.$\end{theorem}
{\bf Proof.}By Holders inequality,we  have
\begin{eqnarray*}
\|h^\alpha_{g}f\|_{\overline A^p(\omega)}^p=\int_{U^n}\omega (1-|z|)\biggl(\int_{U^n}\frac{(1-|\zeta |^2)^\alpha}{|1-\overline z
  \zeta|^{\alpha +2}}|f(\zeta)||g(\zeta )|dm_{2n}(\zeta)\biggr)^pdm_{2n}(z)\leq\\
 \leq C\int_{U^n}\frac{\omega (1-|z|)}{(1-|z|^2)^{\alpha p/q}}
 \int_{U^n}\frac{(1-|\zeta |^2)^{\alpha p}}{|1-\overline z
  \zeta|^{\alpha +2}}|f(\zeta)|^p|g(\zeta )|^pdm_{2n}(\zeta)dm_{2n}(z)\\
  \leq C\int_{U^n}(1-|\xi|^2)^{\alpha p}|f(\zeta)|^p|g(\zeta)|^p\frac{\omega(1-|\zeta|)}{(1-|\zeta|^2)^{\alpha+\alpha p/q}}dm_{2n}(\zeta)= C\|g\|^p_{\infty}\|f\|^p_{A^p(\omega)}
  \end{eqnarray*}
  In the last inequality we have used Lemma 2.
Conversely, let $h^\alpha_g(f)\in \overline A^p(\omega)$ for all $g\in L^\infty (U^n)$. Take the  function $f_r(\zeta )$  and  $g_r(\zeta )$ as follows:
\begin{eqnarray}\label{ef1}f_r(z)=C_r(1-rz)^{-k},
k_j>(\alpha_{\omega_j}+2)/p,\,\, 1\leq j\leq n, 
r=(r_1,...,r_n)\\r_j\in (0,1),\,\,k=(k_1,...,k_n),\,\, C_r=(1-r)^{k-2/p}\omega^{-1/p}(1-r)
\end{eqnarray}
 and
\begin{equation}\label{ef2}g_r(\zeta )=\exp^{-\arg f_r(\zeta )}
\end{equation}

Then we have $\|f_r\|_{ A^p(\omega)}\thicksim const.$ 

We consider the following domain
$$
\widetilde U_j=\{z_j\in U, |\arg z_j|<(1-r_j)/2;\,(4r_j-1)/3<|z_j|<(1+2r_j)/3\}
$$

\begin{equation}\label{ea}
\widetilde U^n=\widetilde U_1\times\ldots \times \widetilde U_n
\end{equation}

 \medskip\noindent the domain  $\widetilde U^n$  defined by (\ref{ea}) andc a polydisk 
 $V^n$  such that $\overline V^n\subset \widetilde U^n$( $\overline V^n$ is the closure of $V^n$)i. Then we get
$$
\|h^\alpha_{g_r}f_r\|_{\overline A^p(\omega)}\geq
\int_{U^n}\omega (1-|z|)\biggl(\int_{ V^n}\frac{(1-|\zeta |)^\alpha}{|1-\overline z\zeta |^{\alpha +2}}|f_r(\zeta
)|dm_{2n}(\zeta)\biggr)^{p}dm_{2n}(z)=I.
$$
 Let 
$$\max_{\substack\zeta\in\overline V^n}|1-\overline z\zeta|=|1-\overline z\widetilde \zeta |,
$$ then 
\begin{eqnarray}
I\geq C_1(\alpha, p, \omega)\frac{(1-r)^{\alpha p-2}}{\omega (1-r)}\int_{U^n}\omega (1-|z|)\biggl(\int_{ V^n}\frac{dm_{2n}(\zeta)}{|1-\overline z\zeta |^{\alpha +2}}\biggr)^pdm_{2n}(z)\geq\nonumber\\
C_1(\alpha, p, \omega)\frac{(1-r)^{(\alpha +2)p-2}}{\omega (1-r)}\int_{U^n}\frac{\omega (1-|z|)}{|1-\overline z\widetilde\zeta|^{(\alpha +2) p}} dm_{2n}(z).\nonumber
\end{eqnarray}
If we assume that $(\alpha_j+2)p\leq\alpha_{\omega_j}+2$ for some $j$, then  for the corresponding integral taking  $\omega_j(t)=t^{\alpha_{\omega_j}}$ we get

$$
\int_{U}\frac{\omega_j (1-|z_j|)}{|1-\overline z\widetilde\zeta|^{(\alpha_j+2) p}} dm_{2}(z_j)\thicksim const,\quad \mathrm{if}\quad 
(\alpha_j+2)p<\alpha_{\omega_j}+2
$$

and
$$
\int_{U}\frac{\omega_j (1-|z_j|)}{|1-\overline z\widetilde\zeta|^{(\alpha_j+2) p}} dm_{2}(z_j)\thicksim \log\frac{1}{1-|\widetilde\zeta_j|},\quad \mathrm{if}\quad (\alpha_j+2)p=\alpha_{\omega_j}+2.
$$

Consequently, 
$$
\frac{(1-r_j)^{(\alpha_j+2) p}}{\omega_j (1-r_j)(1-r_j)^2}\rightarrow \infty,\quad \frac{(1-r_j)^{(\alpha_j+2) p}}{\omega_j (1-r_j)(1-r_j)^2}\log\frac{1}{1-r_j}\rightarrow \infty
$$

$\mathrm{if}\quad r_j\rightarrow 1-0 .$

\hfill
$\square$

Next we consider  the boundedness of the Berezin-typ operators.  
\begin{theorem} Let $p> 1,\,\,f\in A^p(\omega)$ (or $f \in\overline A^p(\omega )) ,\,g\in L^\infty (U^n)$ and let $\alpha_j>(\alpha_{\omega_j}+1)/p-1,\,1\leq j\leq n.$  Then $B^\alpha_g(f)\in L^p(\omega).$  
\end{theorem}
{\bf Proof.} Let $f\in A^p(\omega)$ or $f\in \overline A^p(\omega).$ We will show that $B_\alpha f\in
L^p(\omega ) .$ To this end we estimate the corresponding integral
$$
\int_{U^n}\omega (1-|z|)\biggl((1-|z|^2)^{\alpha
  +2}\int_{U^n}\frac{(1-|\zeta |^2)^\alpha|f(\zeta)|}{|1-z\overline
  \zeta|^{4+2\alpha}}dm_{2n}(\zeta)\biggr)^pdm_{2n}(z)\equiv I
$$
By Holders inequality, we get
\begin{eqnarray*}
\biggl(\int_{U^n}\frac{(1-|\zeta |^2)^\alpha|f(\zeta)|}{|1-z\overline
  \zeta|^{4+2\alpha}}dm_{2n}(\zeta)\biggr)^p\leq\\ \frac{\|g\|_\infty}{(1-|z|)^{(2+\alpha)p/q}}\int_{U^n}\frac{(1-|\zeta|^2)^{\alpha}|f(\zeta)|^p}{|1-\zeta\overline z|^{4+2\alpha}}dm_{2n}(\zeta)
 \end{eqnarray*}
 Then we get
 \begin{eqnarray*}
 I\leq\int_{U^n}\frac{\omega(1-|z|)(1-|z|)^{(\alpha+2)p}}{(1-|z|)^{(\alpha+2)p/q}}\int_{U^n}\frac{(1-|\zeta|^2)^\alpha|f(\zeta)|^p|g(\zeta)|^p}{|1-\zeta\overline z|^{4+2\alpha}}dm_{2n}(\zeta)dm_{2n}(z)\\
= \int_{U^n}|f(\zeta)g(\zeta)|^p(1-|\zeta|^2)^\alpha\int_{U^n}
 \frac{\omega(1-|z|)(1-|z|)^{(\alpha+2)(1-1/q)p}}{|1-\zeta\overline z|^{4+2\alpha}}dm_{2n}(z)dm_{2n}(\zeta)\leq\\
 \int_{U^n}|f(\zeta)g(\zeta)|^p\omega(1-|\zeta|)dm_{2n}{\zeta}\leq \|f\|^p_{A^p(\omega)}\|g\|^p_\infty
 \end{eqnarray*}

We have used Lemma 2 again.

\hfill
$\square$

\end{document}